\documentclass[12pt]{amsart}
\usepackage[nobysame,abbrev,alphabetic]{amsrefs}
\usepackage{amssymb}
\usepackage[all]{xy}
\usepackage{stmaryrd}

\DeclareMathOperator{\Char}{Char}

\DeclareMathOperator{\Gal}{Gal}

\DeclareMathOperator{\id}{id}
\DeclareMathOperator{\Img}{Im}

\DeclareMathOperator{\invlim}{\varprojlim}
\DeclareMathOperator{\Ker}{Ker}

\DeclareMathOperator{\Path}{Path}

\DeclareMathOperator{\Spec}{Spec}

\DeclareFontFamily{U}{wncy}{}
\DeclareFontShape{U}{wncy}{m}{n}{<->wncyr10}{}
\DeclareSymbolFont{mcy}{U}{wncy}{m}{n}
\DeclareMathSymbol{\Sha}{\mathord}{mcy}{"58}
\DeclareMathSymbol{\sha}{\mathord}{mcy}{"78}

\begin{document}

\newtheorem{thm}{Theorem}[section]
\newtheorem{cor}[thm]{Corollary}
\newtheorem{lem}[thm]{Lemma}
\newtheorem{prop}[thm]{Proposition}
\newtheorem{defin}[thm]{Definition}
\newtheorem{exam}[thm]{Example}
\newtheorem{examples}[thm]{Examples}
\newtheorem{rem}[thm]{Remark}
\newtheorem{case}{\sl Case}
\newtheorem{claim}{Claim}
\newtheorem{prt}{Part}
\newtheorem*{mainthm}{Main Theorem}
\newtheorem*{thmA}{Theorem A}
\newtheorem*{transfer theorem}{The Transfer Theorem}

\newtheorem{question}[thm]{Question}
\newtheorem*{notation}{Notation}
\swapnumbers
\newtheorem{rems}[thm]{Remarks}
\newtheorem*{acknowledgment}{Acknowledgment}

\newtheorem{questions}[thm]{Questions}
\numberwithin{equation}{section}

\newcommand{\ab}{\mathrm{ab}}
\newcommand{\Ann}{\mathrm{Ann}}
\newcommand{\Bock}{\mathrm{Bock}}
\newcommand{\dec}{\mathrm{dec}}
\newcommand{\diam}{\mathrm{diam}}
\newcommand{\dirlim}{\varinjlim}
\newcommand{\discup}{\ \ensuremath{\mathaccent\cdot\cup}}
\newcommand{\divis}{\mathrm{div}}
\newcommand{\et}{\mathrm{et}}
\newcommand{\gr}{\mathrm{gr}}
\newcommand{\nek}{,\ldots,}
\newcommand{\ind}{\hbox{ind}}
\newcommand{\Ind}{\mathrm{Ind}}
\newcommand{\inv}{^{-1}}
\newcommand{\isom}{\cong}
\newcommand{\Massey}{\mathrm{Massey}}
\newcommand{\ndiv}{\hbox{$\,\not|\,$}}
\newcommand{\nil}{\mathrm{nil}}
\newcommand{\pietale}{\pi_1^{\mathrm{\acute{e}t}}}
\newcommand{\pr}{\mathrm{pr}}
\newcommand{\sep}{\mathrm{sep}}
\newcommand{\sh}{\mathrm{sh}}
\newcommand{\tagg}{^{''}}
\newcommand{\tensor}{\otimes}
\newcommand{\UT}{\mathbb{UT}}
\newcommand{\alp}{\alpha}
\newcommand{\gam}{\gamma}
\newcommand{\Gam}{\Gamma}
\newcommand{\del}{\delta}
\newcommand{\Del}{\Delta}
\newcommand{\eps}{\epsilon}
\newcommand{\lam}{\lambda}
\newcommand{\Lam}{\Lambda}
\newcommand{\sig}{\sigma}
\newcommand{\Sig}{\Sigma}
\newcommand{\bfA}{\mathbf{A}}
\newcommand{\bfB}{\mathbf{B}}
\newcommand{\bfC}{\mathbf{C}}
\newcommand{\bfF}{\mathbf{F}}
\newcommand{\bfP}{\mathbf{P}}
\newcommand{\bfQ}{\mathbf{Q}}
\newcommand{\bfR}{\mathbf{R}}
\newcommand{\bfS}{\mathbf{S}}
\newcommand{\bfT}{\mathbf{T}}
\newcommand{\bfZ}{\mathbf{Z}}
\newcommand{\dbA}{\mathbb{A}}
\newcommand{\dbC}{\mathbb{C}}
\newcommand{\dbF}{\mathbb{F}}
\newcommand{\dbG}{\mathbb{G}}
\newcommand{\dbI}{\mathbb{I}}
\newcommand{\dbK}{\mathbb{K}}
\newcommand{\dbN}{\mathbb{N}}
\newcommand{\dbP}{\mathbb{P}}
\newcommand{\dbQ}{\mathbb{Q}}
\newcommand{\dbR}{\mathbb{R}}
\newcommand{\dbU}{\mathbb{U}}
\newcommand{\dbV}{\mathbb{V}}
\newcommand{\dbZ}{\mathbb{Z}}
\newcommand{\grf}{\mathfrak{f}}
\newcommand{\gra}{\mathfrak{a}}
\newcommand{\grA}{\mathfrak{A}}
\newcommand{\grB}{\mathfrak{B}}
\newcommand{\grh}{\mathfrak{h}}
\newcommand{\grH}{\mathfrak{H}}
\newcommand{\grI}{\mathfrak{I}}
\newcommand{\grL}{\mathfrak{L}}
\newcommand{\grm}{\mathfrak{m}}
\newcommand{\grp}{\mathfrak{p}}
\newcommand{\grq}{\mathfrak{q}}
\newcommand{\grr}{\mathfrak{r}}
\newcommand{\grR}{\mathfrak{R}}
\newcommand{\grU}{\mathfrak{U}}
\newcommand{\grZ}{\mathfrak{Z}}
\newcommand{\calA}{\mathcal{A}}
\newcommand{\calB}{\mathcal{B}}
\newcommand{\calC}{\mathcal{C}}
\newcommand{\calE}{\mathcal{E}}
\newcommand{\calG}{\mathcal{G}}
\newcommand{\calH}{\mathcal{H}}
\newcommand{\calJ}{\mathcal{J}}
\newcommand{\calK}{\mathcal{K}}
\newcommand{\calL}{\mathcal{L}}
\newcommand{\calR}{\mathcal{R}}
\newcommand{\calW}{\mathcal{W}}
\newcommand{\calU}{\mathcal{U}}
\newcommand{\calV}{\mathcal{V}}
\newcommand{\calZ}{\mathcal{Z}}

\title{Generalized Steinberg Relations}

\author{ Ido Efrat}
\address{Earl Katz Family Chair in Pure Mathematics\\
Department of Mathematics\\
Ben-Gurion University of the Negev\\
P.O.\ Box 653, Be'er-Sheva 84105\\
Israel} \email{efrat@bgu.ac.il}

\thanks{This work was supported by the Israel Science Foundation (grant No.\ 569/21).
\\ \null \ \ \ \ Data sharing not applicable to this article as no datasets were generated or analysed during the current study.
Author declares no conflict of interests.}



\maketitle

\begin{abstract}
\quad
We consider a field $F$ and positive integers $n,m$, such that $m$ is not divisible by $\mathrm{Char}(F)$ and is prime to $n!$.
The absolute Galois group $G_F$ acts on the group $\mathbb{U}_n(\mathbb{Z}/m)$ of all $(n+1)\times(n+1)$ unipotent upper-triangular matrices over $\mathbb{Z}/m$ cyclotomically.
Given $0,1\neq z\in F$ and an arbitrary list $w$ of $n$ Kummer elements $(z)_F$, $(1-z)_F$ in $H^1(G_F,\mu_m)$, we construct in a canonical way a quotient $\mathbb{U}_w$ of $\mathbb{U}_n(\dbZ/m)$ and a cohomology element $\rho^z$ in $H^1(G_F,\mathbb{U}_w)$ whose projection to the superdiagonal is the prescribed list.
This extends results by Wickelgren, and in the case $n=2$ recovers the Steinberg relation in Galois cohomology, proved by Tate.
\end{abstract}
\tableofcontents

\section{Introduction}
Let $F$ be a field with separable closure $F^{\mathrm{sep}}$ and absolute Galois group $G_F=\Gal(F^{\mathrm{sep}}/F)$.
Let $m$ be a positive integer prime to the characteristic $\Char(F)$ of $F$.
We denote the group of $m$th roots of unity in $F^{\mathrm{sep}}$ by $\mu_m$, and write $H^l(G_F,\mu_m^{\tensor l})$ for the $l$th Galois cohomology group with the $l$-times twisted cyclotomic action of $G_F$ on $\mu_m$.
The map of exponentiation by $m$ gives a short exact sequence of discrete $G_F$-modules
\[
1\to\mu_m\to (F^{\mathrm{sep}})^\times\xrightarrow{m}(F^{\mathrm{sep}})^\times\to1.
\]
The corresponding cohomology exact sequence gives rise in particular to the \textsl{Kummer homomorphism}
\[
F^\times\to H^1(G_F,\mu_m), \quad z\mapsto (z)_F.
\]
Its kernel is $(F^\times)^m$, and since, by Hilbert Theorem 90, $H^2(G_F,(F^{\mathrm{sep}})^\times)=0$, it is surjective.

A fundamental fact, proved by Tate \cite{Tate76}, is that the Kummer homomorphism satisfies the  \textsl{Steinberg relation}:
For $0,1\neq z\in F$ one has
\[
(1-z)_F\cup(z)_F=0,
\]
where the cup product is in $H^2(G_F,\mu_m^{\tensor2})$.
However much more is true:
By the celebrated Norm Residue Theorem, proved by Voevodsky and Rost (\cite{Voevodsky03}, \cite{Voevodsky11}), this relation in fact gives rise to \textsl{all} nontrivial relations in $H^l(G_F,\mu_m^{\tensor l})$, $l\geq0$, in the following sense:
Let $K^M_l(F)=(F^\times)^{\tensor l}/ {\rm St}_l(F)$ be the \textsl{$l$th Milnor $K$-group} of $F$, where ${\rm St}_l(F)$ is the subgroup of the tensor power $(F^\times)^{\tensor l}$ generated by all tensors $a_1\tensor\cdots\tensor a_l$ such that $a_i+a_j=1$ for some $1\leq i<j\leq l$.
In view of the Steinberg relations, the Kummer homomorphism induces via the cup product a homomorphism
\[
K^M_l(F)/m\xrightarrow{\sim} H^l(G_F,\mu_m^{\tensor l}).
\]
The Norm Residue Theorem is the surprising fact that this homomorphism is actually \textsl{an isomorphism}.
Thus the Steinberg relations provide a very simple and completely elementary description of the cohomology ring $\bigoplus_{l\geq0} H^l(G_F,\mu_m^{\tensor l})$.
We refer to  \cite{HaesemeyerWeibel19} for a comprehensive introduction to this theorem and its history.

Now the vanishing of cup products can be expressed in more basic Galois-theoretic terms as embedding problems for homomorphisms, and when the Galois action is nontrivial, for $1$-cocycles.
Namely, for $n\geq1$, let $\dbU_n(\dbZ/m)$ be the group of all upper-triangular unipotent matrices of size $(n+1)\times(n+1)$ over the ring $\dbZ/m$.
It carries a natural cyclotomic $G_F$-action, where the copies of $\dbZ/m$ on the $l$-th diagonal above the main one are identified with the $G_F$-module $\mu_m^{\tensor l}$ (see \S\ref{section on equivariant maps}).
We write $\pr_{i,i+1}\colon \dbU_n(\dbZ/m)\to\mu_m$ for the projection on the $(i,i+1)$-entry.
Then, given $\varphi_1,\varphi_2\in H^1(G_F,\mu_m)$, one has $\varphi_1\cup\varphi_2=0$ in $H^2(G_F,\mu_m^{\tensor2})$ if and only if there is a continuous $1$-cocycle $\rho \colon G_F\to\dbU_2(\dbZ/m)$ whose cohomology class $[\rho]$ satisfies $\varphi_i=(\pr_{i,i+1})_*[\rho]$, $i=1,2$ (see Example \ref{example n 2}).
When $\mu_m\subseteq F$, this is a solution of a profinite embedding problem for $G_F$ in the usual sense \cite{FriedJarden08}*{\S22.3}.

Thus the Steinberg relation means that for every $0,1\neq z\in F$ there is a continuous $1$-cocycle $\rho^z\colon G_F\to\dbU_2(\dbZ/m)$ such that $(\pr_{1,2})_*[\rho^z]=(1-z)_F$ and $(\pr_{2,3})_*[\rho^z]=(z)_F$.

In her remarkable work \cite{Wickelgren12b}, Kirsten Wickelgren extends this for $m$ odd as follow:
She proves that, for $n\geq2$ and given $1\leq k\leq n$, one can associate  in a canonical way to every $z\in F^\times$ a continuous $1$-cocycle $\rho^z\colon G_F\to\dbU_n(\dbZ/m)$ such that $(\pr_{i,i+1})_*[\rho^z]$ is $(1-z)_F$ for $i=k$, and is $(z)_F$ for all other values of $1\leq i\leq n$.
Here it is assumed that $\Char(F)=0$ and that $m$ is prime to $n!$.
In subsequent work \cite{Wickelgren17}, she also constructs (in an arbitrary characteristic not dividing $m$) a $1$-cocycle $\rho^z$ such that $(\pr_{i,i+1})_*[\rho^z]$ is $(1-z)_F$ for $i=1$ and $i=n$ and is $(z)_F$ for $2\leq i\leq n-1$.
However this time the canonical continuous $1$-cocycle $\rho^z$ is into a quotient group $\overline\dbU_n(\dbZ/m)$ of $\dbU_n(\dbZ/m)$, which consists of all unipotent upper-triangular $(n+1)\times(n+1)$-matrices over $\dbZ/m$ but with the $(1,n+1)$-entry omitted.
Moreover, an additional construction uses this to obtain a $1$-cocycle into the full group $\dbU_n(\dbZ/m)$.

These facts can be interpreted in the language of external cohomological operations, as descriptions of elements of the $n$-fold Massey products in $H^2(G_F,\mu_m^{\tensor2})$ of the appropriate lists of $(1-z)_F$ and $(z)_F$ -- see \S\ref{section on Massey Products}.

In the current paper we extend these results by Tate and Wickelgren further to \textsl{arbitrary lists} of Kummer elements $(1-z)_F$ and $(z)_F$.
This full generality is achieved at the cost of restricting the target matrix group $\dbU_n(\dbZ/m)$ to a natural quotient, which depends on our specific list.
More specifically, we consider a field $F$ and assume that $m$ is not divisible by $\Char(F)$ and is prime to $n!$.
For notational convenience we encode such lists by words $w=(a_1\cdots a_n)$ in the two-letter alphabet $\{x,y\}$, where $x$ stands for $(z)_F$, and $y$ stands for $(1-z)_F$.
Let $\dbU_w$ be the group obtained from $\dbU_n(\dbZ/m)$ by restricting only to entries $(i,j)$ such that $i\leq j$ and the letters $a_i,a_{i+1}\nek a_{j-1}$ are all $x$ with at most one exception.
Thus these entries include the diagonal, the superdiagonal, and usually larger portions of the upper-right triangle of the matrices (see \S\S\ref{section on Unitriangular Matrices}-\ref{section on Words}, where $\dbU_w$ is denoted $\dbU_{\widetilde T_{w,x}}$).
We prove:

\begin{mainthm}
For every $0,1\neq z\in F$ one can associate in a canonical way a cohomology class $\rho_w^z\in H^1(G,\dbU_w)$ such that, in $H^1(G,\mu_m)$,
\[
(\pr_{i,i+1})_*(\rho_w^z)
=\begin{cases}(z)_F,&\hbox{if }a_i=x,\\
(1-z)_F,&\hbox{if }a_i=y,
\end{cases}
\qquad i=1,2\nek n.
\]
\end{mainthm}

Thus for $w=(yx)$ we have $\dbU_w=\dbU_2(\dbZ/m)$, and we recover Tate's original relation for $m$ odd.
When $w=(x\cdots xyx\cdots x)$ we have $\dbU_w=\dbU_n(\dbZ/m)$, and we recover Wickelgren's  result in \cite{Wickelgren12b}.
Finally, when $w=(yxx\cdots xxy)$ we have $\dbU_w=\overline{\dbU}_n(\dbZ/m)$, recovering the above-mentioned result from \cite{Wickelgren17}.
More complicated words give other groups of partial upper-triangular unipotent matrices.
For instance, when $w=(yxxyy)$ the group $\dbU_w$ consists of all partial matrices of the form
\[
\begin{bmatrix}
1&*&*&*&-&-\\
&1&*&*&*&-\\
&&1&*&*&-\\
&&&1&*&-\\
&&&&1&*\\
&&&&&1
\end{bmatrix},
\]
where $*$ stands for an arbitrary element of $\dbZ/m$, and $-$ for a deleted entry (see Example \ref{concrete example}).

As in \cite{Wickelgren12b} and \cite{Wickelgren17}, we also obtain an analogous result, where the Kummer elements $(z)_F,(1-z)_F$ are replaced by the elements $(z)_F,(-z)_F$ (see Remark \ref{t -t}).
When $n=2$ this recovers the identity $(-z)_F\cup(z)_F=0$, which is a formal consequence of Tate's relation \cite{EfratBook}*{Prop.\ 24.1.2(b)}.

The proof of the Main Theorem combines Wickelgren's method with a new construction of a graded Lie algebra $\gr_\bullet(Z,n)$ related to filtrations of $\dbU_n(\dbZ/m)$.
More specifically, let $\calC$ be the formation of all finite groups of order prime to $\Char(F)$ and to $n!$.
We consider the pro-$\calC$ \'etale fundamental group $S=\pietale((\dbP^1\setminus\{0,1,\infty\})_{\bar F},\overline{01})^\calC$  of the projective line minus three points with the standard tangential base point $\overline{01}$.
It is a free pro-$\calC$ group on two generators $x,y$, and the natural action of $G_F$ on $S$ is well-understood (see \S\ref{section on The Etale Fundamental Group}).
A main step is to construct a continuous homomorphism $\bar\varphi_w\colon S\to\dbU_w$ which is \textsl{$G_F$-equivariant}, that is, it respects the $G_F$-action.
Its composition with the nonabelian Kummer map from the rational base points of $\dbP^1\setminus\{0,1,\infty\}$ to the nonabelian cohomology set $H^1(G_F,S)$ has the desired properties (\S\S\ref{section on The Etale Fundamental Group}-\ref{section on Proof of the Main Theorem}).
It is the graded Lie algebra construction which identifies the precise quotient $\dbU_w$ of $\dbU_n(\dbZ/m)$ for which $\bar\varphi_w$ is $G_F$-equivariant.
This construction (given in  \S\ref{section on Unitriangular Matrices}) is of an elementary combinatorial nature, and associates to any subset $Z$ of $\{1,2\nek n\}$ a filtration of $\dbU_n(\dbZ/m)$ by normal subgroups.
The distribution of letters in the word $w$ gives rise to special sets $Z$ (see \S\ref{section on Words}), and this leads to the definition of $\dbU_w$ and the equivariance of $\bar\varphi_w$.

In \S\ref{section on Massey Products} we generalize the notion of Massey products in profinite group cohomology to a natural context arising from the Lie algebra $\gr_\bullet(Z,n)$.
This extends the interpretation, due to Dwyer \cite{Dwyer75}, of Massey products in $H^2$ in terms of the groups $\dbU_n$ and $\overline\dbU_n$.

\medskip

Massey products for absolute Galois groups of fields were the subject of extensive research in recent years.
Hopkins and Wickelgren \cite{HopkinsWickelgren15} proved the triviality of 3-fold Massey products over local and global fields, and Min\'a\v c and T\^an conjectured that $n$-fold Massey products should be trivial (in an appropriate sense) for every $n\geq3$ and over all fields \cite{MinacTan15}.
So far this was proved for $n=3$  (\cite{Matzri14}, \cite{EfratMatzri16}, \cite{MinacTan16}),  for number fields  \cite{HarpazWittenberg19}, and for fields $F$ such that the maximal pro-$p$ Galois group $G_F(p)$ has $p$-elementary type \cite{Quadrelli22}.
We refer to \cite{HarpazWittenberg19} for more details on the history and the current developments in this direction, as well as for additional references.

In particular, some of these recent works on Massey products in Galois cohomology include results that partly overlap special instances of  \cite{Wickelgren12b}:

In  \cite{MinacTan16}*{\S5}, Min\'a\v c and T\^an consider the case where $m=p$ is a prime number larger than $n$, and $G$ is a profinite group acting trivially on $\dbZ/p$.
Given continuous homomorphisms $\varphi_1,\varphi_2\colon G\to \dbZ/p$ with $\varphi_1\cup\varphi_2=0$ they construct a continuous homomorphism $\rho\colon G\to\dbU_n(\dbZ/m)$ whose  projections to the superdiagonal are all $\varphi_1$, with at most one exception which is $\varphi_2$ appearing at place $1,2,n-1$, or $n$.
Thus, using Tate's relation and under these assumptions, one recovers the main result of \cite{Wickelgren12b}  for the words
\[
(xx\cdots xx),  \  (yxx\cdots xx), \  (xyxx\cdots xx), \  (xx\cdots xyx),  \  (xx\cdots xxy).
\]
It may be interesting to know if our combinatorial techniques can yield in a \textsl{purely algebraic way} homomorphisms (or more generally, $1$-cocycles) in additional situations.

Sharifi \cite{Sharifi07} obtains in an Iwasawa-theoretic context, similar results on the vanishing of Massey products for the situations $(xx\cdots xxy)$,  $(yxx\cdots xx )$.

In the recent work \cite{HayLamSharifiWangWake20}, Hay et al.~introduce a generalized approach to Massey products and related cohomological constructions, which
extend the vanishing results of \cite{Sharifi07}.
Moreover, they show that the vanishing of Massey products in the situation $(xx\cdots xxy)$, when $m=p$ is odd, implies the (known) triviality of all 3-fold Massey products in $H^2(G_F,\dbZ/p)$.
They also raise the interesting question whether  the vanishing of Massey products in the $(xx\cdots xxy)$ situation implies the vanishing of all nonempty $n$-fold Massey products in $H^2(G_F,\dbZ/p)$, $n\geq3$,  for $m=p$ prime such that $\mu_p\subseteq F$ \cite{HayLamSharifiWangWake20}*{\S1.4}.

In a similar spirit, Matzri shows in \cite{Matzri18}, that when $m=p$ is an odd prime, the triviality of nonempty 3-fold Massey products in the $(xyx)$ situation implies the triviality of all nonempty 3-fold Massey products.

\medskip

I thank Marco Antei, Ishai Dan-Cohen and Mathieu Florence for valuable discussions on this subject.
I thank the referees for their careful reading of this paper, and their comments and suggestions which improved the final version.

\section{Unipotent Upper-Triangular Matrices}
\label{section on Unitriangular Matrices}
Let $n\geq2$ be an integer.
We denote
\[
\dbI=\bigl\{(i,j)\ \bigm|\ 1\leq i\leq j\leq n+1\bigr\}, \quad
\Del=\bigl\{(i,i)\ \bigm|\ 1\leq i\leq n+1\bigr\},
\]
and consider these sets as the entries of the upper-right triangle, resp., the diagonal, of an $(n+1)\times(n+1)$-matrix.
We say that a subset $T$ of $\dbI$ is \textsl{convex} if:
\begin{enumerate}
\item[(1)]
$\Del\subseteq T$;
\item[(2)]
For every $(i,j)\in T$ and $(i',j')\in \dbI$ such that $i\leq i'$ and $j'\leq j$ also $(i',j')\in T$.
\end{enumerate}
Note that unions and intersections of convex subsets of $\dbI$ are also convex.

We fix a profinite ring $\calA$.
Let $\dbU=\dbU_n(\calA)$ be the profinite group of all unipotent upper-triangular $(n+1)\times(n+1)$-matrices over $\calA$,
that is, upper-triangular matrices which are $1$ on the main diagonal.
For a convex subset $T$ of $\dbI$ and for $M=(m_{ij})\in\dbU$, we have a partial matrix $(m_{ij})_{(i,j)\in T}$.
Let
\[
\dbU_T=\dbU_{n,T}(\calA)
\]
be the collection of all such partial matrices.
By the convexity, the profinite group structure of $\dbU$ induces a profinite group structure on $\dbU_T$,
namely,
\[
(m_{ij})_{(i,j)\in T}\cdot(n_{ij})_{(i,j)\in T}=\Bigl(\sum_{i\leq k\leq j}m_{ik}n_{kj}\Bigr)_{(i,j)\in T}.
\]
In particular, $\dbU\isom \dbU_\dbI$.
We set
\[
V_T=\Ker(\dbU\to\dbU_T).
\]

For integers $i<j$ we write $[\![i,j-1]\!]$ for the set of integers $i\leq k<j$.
Let $Z$ be a subset of $\{1,2\nek n\}$, and $t\geq1$ an integer.
We define
\[
T(Z,t)=\bigl\{(i,j)\in \dbI\ \bigm|\ |[\![i,j-1]\!]\cap Z|<t\bigr\}.
\]
Clearly, $T(Z,t)$ is convex.

\begin{lem}
\label{increasing sequence}
For $t_0=|Z|$ one has
\[
T(Z,1)\subsetneq T(Z,2)\subsetneq\cdots\subsetneq T(Z,t_0)\subsetneq \dbI=T(Z,t_0+1)=T(Z,t_0+2)=\cdots.
\]
\end{lem}
\begin{proof}
It is immediately clear that the sequence $T(Z,t)$, $t=1,2\nek$ is weakly increasing.
One has $T(Z,t)=\dbI$ if and only if $(1,n+1)\in T(Z,t)$, which means that $t_0<t$.

When $1\leq t\leq t_0$, we take $j$ such that $j-1$ is the $t$-th smallest number in $Z$.
Then  $(1,j)\in T(Z,t+1)\setminus T(Z,t)$.
\end{proof}

Consequently, we have
\begin{equation}
\label{decreasing sequence of the V subgroups}
V_{T(Z,1)}>V_{T(Z,2)}>\cdots>V_{T(Z,t_0)}>\{I\}=V_{T(Z,t_0+1)}=V_{T(Z,t_0+2)}=\cdots,
\end{equation}
where $I$ is the identity matrix.

\begin{prop}
\label{commutativity of product}
Let $Z$ be a subset of $\{1,2\nek n\}$ and let $t,s\geq1$ be integers.
For $M\in V_{T(Z,t)}$ and $N\in V_{T(Z,s)}$, the commutator $[M,N]$ belongs to $V_{T(Z,t+s)}$.
\end{prop}
\begin{proof}
Write $M=(m_{ij})$ and $N=(n_{ij})$.
We need to show that the images of $MN$ and $NM$ in $\dbU_{T(Z,t+s)}$ coincide, i.e., $(MN)_{ij}=(NM)_{ij}$ for every $(i,j)\in \dbI$ such that $|[\![i,j-1]\!|\cap Z|<t+s$.

Consider $k$ such that $i<k<j$.

If $|[\![i,k-1]\!]\cap Z|<t$, then $(i,k)\in T(Z,t)\setminus\Del$, so $m_{ik}=0$.
If not, then $|[\![k,j-1]\!]\cap Z|<s$, so $(k,j)\in T(Z,s)\setminus\Del$, and we have $n_{kj}=0$.
In both cases $m_{ik}n_{kj}=0$.

Similarly, if $|[\![k,j-1]\!]\cap Z|<t$, then $(k,j)\in T(Z,t)\setminus\Del$, so $m_{kj}=0$.
If not, then $|[\![i,k-1]\!]\cap Z|<s$, so $(i,k)\in T(Z,s)\setminus\Del$, and we have $n_{ik}=0$.
In both cases $n_{ik}m_{kj}=0$.

Consequently,
\[
(MN)_{ij}=\sum_{k=i}^jm_{ik}n_{kj}=m_{ij}+n_{ij}=\sum_{k=i}^jn_{ik}m_{kj}=(NM)_{ij}.
\qedhere
\]
\end{proof}

\begin{cor}
If $|Z|<2t$, then $V_{T(Z,t)}$ is abelian.
\end{cor}

We recall that, for a decreasing sequence of subgroups $V_t$, $t=1,2\nek$ of a group $V$, such that $V_1=V$ and $[V_t,V_s]\leq V_{t+s}$ for every $t,s\geq1$, the quotients $V_t/V_{t+1}$ are abelian, and the commutator map induces on $\bigoplus_{t\geq1}V_t/V_{t+1}$ a graded Lie algebra structure \cite{SerreLie}*{Part I, Ch.\ 2, Prop.\ 2.3 and Prop.\ 3.1}.
In view of Proposition \ref{commutativity of product}, we deduce in particular that the quotient groups $V_{T(Z,t)}/V_{T(Z,t+1)}$ are abelian and the commutator map endows the $\dbZ$-module
\[
\gr_\bullet (Z,n)=\bigoplus_{t\geq1}V_{T(Z,t)}/V_{T(Z,t+1)}
\]
with a structure of a graded Lie algebra.
By (\ref{decreasing sequence of the V subgroups}), the nonzero components of $\gr_\bullet(Z,n)$ are exactly in degrees $1\leq t\leq|Z|$.

This construction is functorial, in the sense that for subsets $Z'\subseteq Z$ of $\{1,2\nek n\}$ one has $V_{T(Z',t)}\leq V_{T(Z,t)}$ for every $t\geq1$, and we obtain a natural graded  Lie algebra homomorphism
$\gr_\bullet(Z',n)\to\gr_\bullet(Z,n)$.

\section{Words}
\label{section on Words}
Let $X$  be a set, considered as an alphabet.
Let $w=(a_1\cdots a_n)$ be a word of length $n$ in the alphabet $X$.
For $x\in X$ let
\[
Z_{w,x}=\{1\leq i\leq n\ |\ a_i\neq x\}.
\]
We set
\[
T_{w,x}=T(Z_{w,x},1)=\bigl\{(i,j)\in \dbI\ \bigm|\ a_i=\cdots=a_{j-1}=x\bigr\}.
\]
For $x,y\in X$, $x\neq y$, one has $T_{w,x}\cap T_{w,y}=\Del$.
The union
\begin{equation}
\label{Tw}
T_w=\bigcup_{x\in X}T_{w,x}=\{(i,j)\in\dbI\ |\ a_i=\cdots=a_{j-1}\},
\end{equation}
is also convex in $\dbI$.

We further set
\[
\widetilde T_{w,x}=T(Z_{w,x},2).
\]
It consists of all $(i,j)\in\dbI$ such that $a_i\nek a_{j-1}$ are $x$ with at most one exception.
Alternatively, let $l_1<\cdots< l_r$ be the places of letters $\neq x$ in $w$.
Then
\[
\widetilde T_{w,x}=\dbI\setminus\bigcup_{s=1}^{r-1}\{(i,j)\in\dbI\  |\   i\leq l_s,\ l_{s+1}<j\},
\]
which is $\dbI$ with $r-1$ upper-right rectangles removed.
We note that $\widetilde T_{w,x}$ contains both the diagonal $\Del$ and the superdiagonal $(i,i+1)$, $i=1,2\nek n$.

\begin{exam}
\rm
\begin{enumerate}
\item[]
\item[(1)]
One has  $\widetilde T_{w,x}=\dbI$ if and only if $w$ contains at most one letter $\neq x$.
\item[(2)]
One has $\widetilde T_{w,x}=\dbI\setminus\{(1,n+1)\}$ if and only if $w=(yxx\cdots xxz)$ where $y,z$ are letters $\neq x$ (possibly $y=z$).
\end{enumerate}
\end{exam}

\begin{exam}
\label{concrete example}
\rm
In the two-letter alphabet $X=\{x,y\}$ let $w=(yxxyy)$.
Then $n=5$, $Z_{w,x}=\{1,4,5\}$, $Z_{w,y}=\{2,3\}$, and the above convex subsets of $\dbI$ can be visualized as
\[
T_{w,x}=\begin{bmatrix}
*&&&&&\\
&*&*&*&&\\
&&*&*&&\\
&&&*&&\\
&&&&*&\\
&&&&&*
\end{bmatrix}, \quad
T_{w,y}=\begin{bmatrix}
*&*&&&&\\
&*&&&&\\
&&*&&&\\
&&&*&*&*\\
&&&&*&*\\
&&&&&*
\end{bmatrix}
\]
\[
\
\tilde T_{w,x}=\begin{bmatrix}
*&*&*&*&&\\
&*&*&*&*&\\
&&*&*&*&\\
&&&*&*&\\
&&&&*&*\\
&&&&&*
\end{bmatrix}, \quad
\tilde T_{w,y}=\begin{bmatrix}
*&*&*&&&\\
&*&*&&&\\
&&*&*&*&*\\
&&&*&*&*\\
&&&&*&*\\
&&&&&*
\end{bmatrix}.
\]
\end{exam}

\medskip

Next let $T$ be an arbitrary convex subset of $\dbI$.
To a matrix $M=(m_{ij})\in\dbU$ we associate a matrix $R_T(M)\in\dbU$ defined by
\[
R_T(M)_{ij}=\begin{cases}m_{ij},& \hbox{if } (i,j)\in T,\\
                                 0,&  \hbox{otherwise.}\end{cases}
\]

\begin{rem}
\label{Im R contained in V}
\rm
If $T,T'$ are convex subsets of $\dbI$, then $\Img(R_{T'})\cap V_{T\cap T'}\subseteq V_T$.
\end{rem}

\begin{lem}
\label{RT a homomorphism}
Let $w=(a_1\cdots a_n)$ be a word of length $n$ in $X$,  let $x\in X$, and set $T=T_{w,x}$.
Then:
\begin{enumerate}
\item[(a)]
$R_T\colon \dbU\to\dbU$ is a homomorphism.
\item[(b)]
$\dbU=V_T\rtimes\dbU_T$.
\end{enumerate}
\end{lem}
\begin{proof}
(a) \quad
For an arbitrary convex subset $T$ of $\dbU$, the projection $\dbU\to\dbU_T$, $M\mapsto\bar M$, is an epimorphism.
In particular, for $M,N\in\dbU$, the matrices $R_T(MN)$ and $R_T(M)R_T(N)$ coincide on $T$.

Now take $T=T_{w,x}$ and $(i,j)\in\dbI\setminus T$.
Thus $a_i\nek a_{j-1}$ are not all $x$.
Hence for every $i\leq k\leq j$, the letters $a_i\nek a_{k-1}$ are not all $x$, or else $a_k\nek a_{j-1}$ are not all $x$.
In the first case, $(i,k)\not\in T$, so $R_T(M)_{ik}=0$, and in the second case $(k,j)\not\in T$, so $R_T(N)_{kj}=0$.
We conclude that
\[
(R_T(M)R_T(N))_{ij}=\sum_{i=k}^j R_T(M)_{ik}R_T(N)_{kj}=0=R_T(MN)_{ij}.
\]
(b) \quad
By (a), the epimorphism $\dbU\to\dbU_T$, $M\mapsto\bar M$, has a homomorphic section $\bar M\mapsto R_T(M)$.
\end{proof}

\section{Equivariant maps}
\label{section on equivariant maps}
Let $\Sig$ be a set of prime numbers.
A \textsl{$\Sig$-integer} is a positive integer whose prime divisors are in $\Sig$.
Let $\calC=\calC(\Sig)$ be the class of all finite groups of $\Sig$-integer order.
Let $\hat\dbZ(\calC)$ be the free pro-$\calC$ cyclic ring, i.e., $\hat\dbZ(\calC)=\invlim\dbZ/m$, where $m$ ranges over the $\Sig$-integers and with the divisibility partial order.

Let $S=\invlim S_i$ be a pro-$\calC$ group, where $S_i\in \calC$.
We recall that there is a \textsl{pro-$\calC$ exponentiation map} $S\times\hat\dbZ(\calC)^+\to S$, $(\sig,\alp)\mapsto\sig^\alp$, linear in $\hat\dbZ(\calC)^+$, which is defined as follows:
Write $\sig=(\sig_i)_i$ with $\sig_i\in S_i$, and $\alp=(a_m\hbox{ (mod }m))_m$, where the integers $a_m$ satisfy $a_m\equiv a_{m'}\pmod m$ for $\Sig$-integers $m|m'$.
For a given $i$, take $m=|S_i|$.
Then  $\sig_i^{a_m}=\sig_i^{a_{m'}}$ for every $m'$ with $m|m'$, so the sequence $(\sig_i^{a_m})_m$ in $S_i$ stabilizes, and we denote its limit by $\sig_i^\alp$.
These limits are compatible under the inverse system epimorphisms, which yields a well-defined element $\sig^\alp=(\sig_i^\alp)_i$ of $S$.

This construction is functorial, in the sense that for sets $\Sig\supseteq \Sig'$ of prime numbers, for a pro-$\calC(\Sig)$ (resp., pro-$\calC(\Sig')$) group $S$ (resp., $S'$), and for a continuous homomorphism $S\to S'$,
there is a commutative diagram
\begin{equation}
\label{functoriality of profinite powers}
\xymatrix{
S\ar[d]&*-<3pc>{\times}& \hat \dbZ(\calC(\Sig))^+\ar[d]\ar[r]& S\ar[d]\\
S'&*-<3pc>{\times}& \hat \dbZ(\calC(\Sig'))^+\ar[r]& S'.
}
\end{equation}

From now on we fix a pro-$\calC$ cyclic ring $\calA$, where $\calC=\calC(\Sig)$.
Thus there is a ring epimorphism $\hat\dbZ(\calC)\to\calA$.
Let $G$ be a profinite group, and let $\chi\colon G\to\hat\dbZ(\calC)^\times$ be a continuous homomorphism.
Then $G$ acts on $\hat\dbZ(\calC)^+$ and $\calA^+$ continuously and in a compatible way by ${}^g\alp=\chi(g)\alp$.

\begin{lem}
\label{nu is equivariant}
Let $S$ be a pro-$\calC$ group and $X$ a set of (profinite) generators of $S$.
Suppose that $G$ acts continuously on $S$, and ${}^gx$ and $x^{\chi(g)}$ are conjugate in $S$ for every $x\in X$ and $g\in G$.
Then every continuous homomorphism $\nu\colon S\to\calA^+$ is $G$-equivariant.
\end{lem}
\begin{proof}
For every $x\in X$ and $g\in G$ we have by commutativity and (\ref{functoriality of profinite powers}),
\[
\nu({}^gx)= \nu(x^{\chi(g)})=\chi(g)\nu(x)={}^g\nu(x).
\]
Now use the fact that $X$ generates $S$ as a profinite group, and the continuity of the $G$-action on $S$.
\end{proof}

Let $\UT_n(\calA)$ be the ring of all upper-triangular $(n+1)\times(n+1)$-matrices $M$ over $\calA$.
Let $\UT^0_n(\calA)$ be the additive group of all strictly upper-triangular matrices in $\UT_n(\calA)$.
Then $G$ acts on $\UT_n(\calA)$ by
\[
{}^g(m_{ij})=(\chi(g)^{j-i}m_{ij})_{(i,j)\in\dbI}
\]
for $g\in G$ and $M=(m_{ij})\in\UT_n(\calA)$.
In particular, this gives a continuous $G$ action on the pro-$\calC$ groups $\UT^0_n(\calA)$ and $\dbU=\dbU_n(\calA)$, and more generally, on $\dbU_T=\dbU_{n,T}(\calA)$ for every convex subset $T$ of $\dbI$.

The projection homomorphisms $\pr_{i,i+1}\colon\dbU\to\calA^+$ are $G$-equivariant for $1\leq i\leq n$.

For a convex subset $T$ of $\dbI$ and for every $g\in G$ and $M\in \dbU$ we have
\begin{equation}
\label{RT equivariant}
{}^gR_T(M)=R_T({}^gM).
\end{equation}

There is a $G$-equivariant group homomorphism
\begin{equation}
\label{first map}
\calA^+\to \UT^0_n(\calA), \qquad
\alp\mapsto\begin{bmatrix}
0&\alp&0&\cdots&0\\
0&0&\alp&\cdots&0\\
\vdots&\vdots&\vdots&\ddots&\vdots\\
0&0&0&\cdots&\alp\\
0&0&0&\cdots&0\end{bmatrix},
\end{equation}
and matrices of the latter form commute in multiplication.

Assuming that $n!$ is invertible in $\calA$,  there is a well-defined exponential map
\[
\exp\colon\UT^0_n(\calA)\to\dbU, \quad  M\mapsto\sum_{k=0}^n M^k/k!
\]
(note that $M^k=0$ for $k>n$).
It is clearly $G$-equivariant.
It is a special case of the Baker--Campbell--Hausdorff formula \cite{SerreLie}*{Ch.\ IV, \S7} that $\exp(M+N)=\exp(M)\exp(N)$ when $M,N$ commute.
By composing (\ref{first map}) with $\exp$ we obtain a $G$-equivariant continuous homomorphism
\[
\xi\colon \calA^+\to\dbU, \quad
\xi(\alp)=\Biggl(\dfrac{\alp^{j-i}}{(j-i)!}\Biggr)_{(i,j)\in\dbI}.
\]
The fact that $\xi$ is a homomorphism can also be verified directly, using the binomial formula.
Moreover, we note that the restriction of $\xi$ to the subgroup $n!\calA^+$ of $\calA^+$ is well defined even without the invertibility assumption.

\begin{rem}
\rm
In addition to \cite{Wickelgren12b}*{\S2.6} and \cite{Wickelgren17}*{\S2.7}, variants of this construction also appeared in  \cite{MinacTan16}*{\S5}, and \cite{HayLamSharifiWangWake20}*{\S4.2}.
\end{rem}

We deduce from (\ref{RT equivariant}) that for every $\alp\in n!\calA^+$ and $g\in G$,
\begin{equation}
\label{RT xi equivariant}
(R_T\circ\xi)(\chi(g)\alp)={}^g((R_T\circ\xi)(\alp)).
\end{equation}

When $T=T_{w,x}$ for a word $w$ of length $n$ in the alphabet $X$ and for $x\in X$, the map $R_T\circ\xi\colon n!\calA^+\to\dbU$ is a homomorphism of pro-$\calC$ groups (see Lemma \ref{RT a homomorphism}(a)).
By the commutativity of (\ref{functoriality of profinite powers}),
\begin{equation}
\label{RT commutes with chi}
(R_T\circ\xi)(\alp)^{\chi(g)}=(R_T\circ\xi)(\chi(g)\alp)
\end{equation}
for every $\alp\in n!\calA^+$ and $g\in G$.

\section{The homomorphism $\varphi_w$}
\label{section on the Homomorphism varphi w}
Let $\calC$, $\calA$, $G$, and $\chi$ be as in the previous section.
Let $S$ be the free pro-$\calC$ group on a finite basis $X$ \cite{FriedJarden08}*{\S17.4}.
Thus $X$ is a subset of $S$, which generates it as a pro-$\calC$ group, and every map $\gam\colon X\to K$, where $K$ is a pro-$\calC$ group, extends uniquely to a continuous homomorphism $S\to K$.

We fix arbitrary elements
\begin{equation}
\label{alp x}
\alp_x\in n!\calA^+, \quad x\in X.
\end{equation}
We consider $X$ also as an alphabet.
For a word $w=(a_1\cdots a_n)$  in $X$, we define a continuous homomorphism
\[
\varphi_w\colon S\to\dbU, \qquad x\mapsto R_{T_{w,x}}(\xi(\alp_x)) \quad \hbox{     for   } x\in X.
\]

Then, for $g\in G$ and $x\in X$,  (\ref{RT commutes with chi})  and (\ref{RT xi equivariant}) imply that
\begin{equation}
\label{exponentiation equals action}
\varphi_w(x)^{\chi(g)}={}^g\varphi_w(x).
\end{equation}

\begin{rem}
\rm
\cite{Wickelgren12b} and \cite{Wickelgren17} work with $\alp_x=1$ for every $x\in X$, which requires the assumption that $n!$ is invertible in $\calA$.
We will make this assumption only at a later stage, in \S\ref{section on Proof of the Main Theorem}.
\end{rem}

\begin{exam}
\rm
In the setup of Example \ref{concrete example}, the homomorphism $\varphi_w$ is given by
\[
\
x\mapsto\begin{bmatrix}
1&0&0&0&0&0\\
&1&\alp_x&\alp_x^2/2&0&0\\
&&1&\alp_x&0&0\\
&&&1&0&0\\
&&&&1&0\\
&&&&&1
\end{bmatrix}, \quad
y\mapsto\begin{bmatrix}
1&\alp_y&0&0&0&0\\
&1&0&0&0&0\\
&&1&0&0&0\\
&&&1&\alp_y&\alp_y^2/2\\
&&&&1&\alp_y\\
&&&&&1
\end{bmatrix}.
\]
\end{exam}

Let $[S,S]$ be the (closed) commutator subgroup of $S$.
The next fact is essentially \cite{Wickelgren12b}*{Lemma 2.11}.
Here $T_w$ is as in (\ref{Tw}).

\begin{lem}
\label{image of commutator subgroup}
One has $\varphi_w([S,S])\subseteq V_{T_w}$.
\end{lem}
\begin{proof}
As $T_w=\bigcup_{x\in X}T_{w,x}$, we need to show that $\varphi_w([S,S])\subseteq V_{T_{w,x}}$ for every $x\in X$.

To this end, let $N$ be the normal closure  of the set of commutators $[y,z]$, where $y,z\in X$ and $y\neq x,z$.
Thus the cosets of the elements of $X$ commute in $S/N$.
Since $X$ generates $S$, this implies that $S/N$ is abelian, that is, $[S,S]\leq N$.
The opposite inclusion is trivial, so $[S,S]=N$.

Now for every $y\in X$, $y\neq x$, we have $T_{w,x}\cap T_{w,y}=\Del$, so $V_{T_{w,x}\cap T_{w,y}}=\dbU$.
By Remark \ref{Im R contained in V}, $\varphi_w(y)\in V_{T_{w,x}}$.
The normality of $V_{T_{w,x}}$ in $\dbU$ therefore implies that $\varphi_w([y,z])\in V_{T_{w,x}}$ for every $z\in X$.
By the previous observation, this implies that $\varphi_w([S,S])\leq V_{T_{w,x}}$.
\end{proof}

\begin{prop}
\label{G equivariance of maps}
Suppose that $G$ acts on $S$ continuously.
Let $w$ be a word of length $n$ in $X$, and let $x\in X$ and $g\in G$.
\begin{enumerate}
\item[(a)]
If ${}^gx=x^{\chi(g)}$, then $\varphi_w({}^gx)={}^g\varphi_w(x)$.
\item[(b)]
Suppose that ${}^gx=f\inv x^{\chi(g)}f$ for some $f\in [S,S]$.
Then $\varphi_w({}^gx)$ and ${}^g\varphi_w(x)$ coincide on $\widetilde T_{w,y}=T(Z_{w,y},2)$ for every $y\in X$ with $y\neq x$.
\end{enumerate}
\end{prop}
\begin{proof}
(a) \quad
Apply (\ref{functoriality of profinite powers}) and (\ref{exponentiation equals action}).

\medskip

(b) \quad
By (\ref{functoriality of profinite powers}),
$\varphi_w({}^gx)=\varphi_w(f)\inv\varphi_w(x)^{\chi(g)}\varphi_w(f)$.

As $\varphi_w(x)\in \Img(R_{T_{w,x}})$ and $T_{w,x}\cap T_{w,y}=\Del$, Remark \ref{Im R contained in V} provides as before $\varphi_w(x)\in V_{T_{w,y}}=V_{T(Z_{w,y},1)}$.
By Lemma \ref{image of commutator subgroup},  $\varphi_w(f)\in V_{T_w}\leq V_{T_{w,y}}$.
Therefore Proposition \ref{commutativity of product} (with $t=s=1$) implies that $[\varphi_w(x),\varphi_w(f)]\in V_{T(Z_{w,y},2)}$.
Hence the projections to $\dbU_{T(Z_{w,y},2)}$ of $\varphi_w({}^gx)$ and $\varphi_w(x)^{\chi(g)}$ coincide, that is, they coincide on $T(Z_{w,y},2)$.
It remains to apply (\ref{exponentiation equals action}).
\end{proof}

\section{The {\'E}tale Fundamental Group}
\label{section on The Etale Fundamental Group}
In this section we work in profinite non-abelian cohomology, defined similarly to discrete non-abelian cohomology \cite{SerreCG}*{I, \S5}.
Specifically, let $G$ be a profinite group which acts continuously on a profinite group $S$.
A \textsl{$1$-cocycle} $f\colon G\to S$ is a continuous map such that $f(gg')=f(g)\cdot {}^gf(g')$ for every $g,g'\in G$.
We write $Z^1(G,S)$ for the set of all such maps.
Two such $1$-cocycles $f,f'$ are \textsl{cohomologous} if there is $s\in S$ such that $f'(g)=s\inv \cdot f(g)\cdot {}^gs$ for every $g\in G$.
This is clearly an equivalence relation, and the cohomology set $H^1(G,S)$ consists of all equivalence (cohomology)  classes $[f]$.
If $S_1,S_2$ are profinite groups upon which $G$ acts continuously, and $\psi\colon S_1\to S_2$ is a continuous $G$-equivariant homomorphism, then there is a natural map $\psi_*\colon Z^1(G,S_1)\to Z^1(G,S_2)$, $f\mapsto \psi\circ f$, which induces a map $\psi_*\colon H^1(G,S_1)\to H^1(G,S_2)$ on the cohomology.

Now let $F$ be a field with absolute Galois group $G=G_F$.
Let $\Sig$ be a set of prime numbers not containing $\Char(F)$, and let $\calC=\calC(\Sig)$, $\hat\dbZ(\calC)$, $\calA$ be as in \S\ref{section on equivariant maps}.
Let $\chi\colon G\to\hat\dbZ(\calC)^\times$ be the \textsl{cyclotomic character}, i.e.,
${}^g\zeta=\zeta^{\chi(g)}$ for every $g\in G$ and every root of unity $\zeta$ in $\bar F$ whose order is a $\Sig$-integer.

The next discussion follows closely \cite{Wickelgren12a}*{\S12.2.1}, \cite{Wickelgren12b},  and \cite{Wickelgren17}, where one can find more details.
In particular, we refer to these papers, as well as to \cite{Deligne89}*{\S15},  \cite{Nakamura99}, for the notion of  \textsl{tangential base points}.
Note that in  \cite{Wickelgren12a}, \cite{Wickelgren12b} it is assumed that $\Char(F)=0$, however in \cite{Wickelgren17} it is explained how the theory extends to arbitrary characteristics, as long as one assumes (as we do here) that the characteristic is not a prime in $\Sig$.

Let $C$ be a smooth geometrically connected curve over $F$, and write $C_{F^\sep}$ for its extension by scalars to the separable closure $F^\sep$ of $F$.
We write $C^{\rm bp}(F)$ for the set of all \textsl{rational base points} of $C$, that is, either $F$-rational points or tangential base points.
Each such point determines canonically a geometric point in $C_{F^\sep}$, whence a fiber functor with a natural $G$-action.
For $b,b'\in C^{\rm bp}(F)$ let $\Path(b,b')$ be the set of natural transformations (considered as ``paths") from the fiber functor of $b$ to the fiber functor of $b'$.
It has a natural structure of a profinite set.

Now fix $b\in C^{\rm bp}(F)$.
The \textsl{\'etale fundamental group} of $C_{F^\sep}$ with base point $b$ is
\[
\pietale(C_{F^\sep},b)=\Path(b,b).
\]
Let $\pietale(C_{F^\sep},b)^\calC$ be its maximal pro-$\calC$ quotient.
It carries a natural continuous $G$-action.

Given $z\in C^{\mathrm{bp}}(F)$ and $\wp\in \Path(b,z)$, the map
\[
\kappa_{b,\wp}(z)\colon G\to \pietale(C_{F^\sep},b)^\calC, \quad g\mapsto \wp\inv\cdot{}^g\wp
\]
is a continuous $1$-cocycle.
If $\wp'\in\Path(b,z)$ is another path, then
\[
(\wp')\inv\cdot{}^g(\wp') =(\wp\inv\wp')\inv\cdot \wp\inv\cdot{}^g\wp\cdot {}^g(\wp\inv\wp'),
\]
whence $[\kappa_{b,\wp'}(z)]=[\kappa_{b,\wp}(z)]$.
One may therefore define the \textsl{nonabelian Kummer map}
\[
\kappa_b\colon C^{\rm bp}(F)\to H^1(G,\pietale(C_{F^\sep},b)^\calC), \quad z\mapsto [\kappa_{b,\wp}(z)],
\]
with $\wp\in\Path(b,z)$ arbitrary.

Now take
\[
C=\dbP^1\setminus\{0,1,\infty\}=\Spec F\Bigl[t,\frac1t,\frac1{1-t}\Bigr]
\]
 and the standard tangential base point $b=\overline{01}=0+1\eps$.
Let
\[
S=\pietale((\dbP^1\setminus\{0,1,\infty\})_{F^\sep},\overline{01})^\calC.
\]
Then $S$ is a free pro-$\calC$ group on two generators $x,y$, and the Galois action of $G$ on $S$ is given as follows  \cite{Ihara91}:
There is a continuous $1$-cocycle $f\colon G\to [S,S]$ such that for every $g\in G$,
\[
{}^gx=x^{\chi(g)}, \quad  {}^gy=f(g)\inv y^{\chi(g)}f(g).
\]

Define $\nu_x\colon S\to\calA^+$ (resp., $\nu_y\colon S\to\calA^+$) to be the unique continuous homomorphism such that $\nu_x(x)=\alp_x$ and $\nu_x(y)=0$ (resp., $\nu_y(x)=0$ and $\nu_y(y)=\alp_y$), with $\alp_x,\alp_y\in n!\calA^+$ as in (\ref{alp x}).
By Lemma \ref{nu is equivariant}, $\nu_x$ and $\nu_y$ are $G$-equivariant.
We further observe that
\begin{equation}
\label{pr composed with varphi}
\pr_{i,i+1}\circ\varphi_w=\alp_{a_i}\nu_{a_i}, \quad i=1,2\nek n.
\end{equation}

\section{Proof of the Main Theorem}
\label{section on Proof of the Main Theorem}
In this section we assume that
\begin{center}
\textsl{$n!$ is invertible in the ring $\calA$.}
\end{center}
Note that this holds if the primes $\leq n$ are not contained in $\Sig$.
We may therefore take in (\ref{alp x}) $\alp_x=1$ for all $x\in X$.

We consider $X=\{x,y\}$ both as a basis of the free pro-$\calC$ group $S$ as above and as an alphabet.
A main arithmetical fact, shown in  \cite{Wickelgren12b}*{Lemma 3.16} (see also \cite{Wickelgren12a}*{\S12.2.2}, as well as the proof of \cite{Wickelgren17}*{Cor.\ 3.10}) is that then, for every $0,1\neq z\in F$, one has in $H^1(G,\calA^+)$:
\begin{equation}
\label{values of Kummer maps}
(\nu_x)_*(\kappa_b(z))=(z)_F, \qquad (\nu_y)_*(\kappa_b(z))=(1-z)_F.
\end{equation}

Let $w=(a_1\cdots a_n)$ be a word of length $n$ in $X$, and recall from \S\ref{section on Words} that $\widetilde T_{w,x}=T(Z_{w,x},2)$.
Let $\bar\varphi_w\colon S\to\dbU_{\widetilde T_{w,x}}$ be the continuous homomorphism induced from $\varphi_w\colon S\to\dbU$, by composing it with the projection $\dbU\to \dbU_{\widetilde T_{w,x}}$.
It follows from  Proposition \ref{G equivariance of maps} (where in (b) the roles of $x,y$ are interchanged) that $\bar\varphi_w$ is $G$-equivariant.

\begin{thm}
\label{the main theorem}
To every $0,1\neq z\in F$ one can associate in a canonical way a cohomology class $\rho_w^z\in H^1(G,\dbU_{\widetilde T_{w,x}})$ such that in $H^1(G,\calA^+)$,
\[
(\pr_{i,i+1})_*(\rho_w^z)
=\begin{cases}(z)_F,&\hbox{if }a_i=x,\\
(1-z)_F,&\hbox{if }a_i=y,
\end{cases}
\qquad i=1,2\nek n.
\]
\end{thm}
\begin{proof}
Given $0,1\neq z\in F$, we choose $\wp\in\Path(b,z)$.
Then $\kappa_{b,\wp}(z)\in Z^1(G,S)$.
Set $\rho_{w,\wp}^z=\bar\varphi_w\circ(\kappa_{b,\wp}(z))$.
The $G$-equivariance of $\bar\varphi_w$ implies that $\rho_{w,\wp}^z\in Z^1(G,\dbU_{\widetilde T_{w,x}})$.
Therefore
\[
\rho_w^z:=[\rho_{w,\wp}^z]=(\bar\varphi_w)_*(\kappa_b(z))\in H^1(G,\dbU_{\widetilde T_{w,x}}),
\]
and this cohomology class is independent of the choice of $\wp$.

Since $\widetilde T_{w,x}$ contains the superdiagonal, and by (\ref{pr composed with varphi}), $\pr_{i,i+1}\circ\bar\varphi_w=\nu_{a_i}$, $i=1,2\nek n$.
Using (\ref{values of Kummer maps}) we compute:
\[
\begin{split}
(\pr_{i,i+1})_*(\rho_w^z)&=
[\pr_{i,i+1}\circ \bar\varphi_w\circ(\kappa_{b,\wp}(z))]=
[\nu_{a_i}\circ(\kappa_{b,\wp}(z))] \\
&=(\nu_{a_i})_*(\kappa_b(z))
=\begin{cases}(z)_F,&\hbox{if }a_i=x,\\
(1-z)_F,&\hbox{if }a_i=y.
\end{cases}
\end{split}
\]
\end{proof}

The Main Theorem stated in the Introduction is the case where $\Sig$ consists of all prime numbers different from $\Char(F)$ and larger than $n$, $m$ is a $\Sig$-integer, and $\calA=\dbZ/m$ with the cyclotomic action (so $\calA^+=\mu_m$ as a Galois module).
 \cite{Wickelgren12b}, \cite{Wickelgren17} work with $\calA=\hat\dbZ(\Sig)$ where $\Sig$ is a set of primes larger than $n$ and different from $\Char(F)$.

\begin{rem}
\label{automorphism of U T}
\rm
For every convex subset $T$ of $\dbI$, the map $(m_{ij})_{(i,j)\in T}\mapsto((-1)^{j-i}m_{ij})_{(i,j)\in T}$ is an automorphism of $\dbU_T$.
Hence we could have taken in Theorem \ref{the main theorem} the Kummer elements to be $(z\inv)_F, ((1-z)\inv)_F$ instead of  $(z)_F, (1-z)_F$, and this is indeed how the result is stated in \cite{Wickelgren12b}.
\end{rem}

\begin{rem}
\label{t -t}
\rm
The ring automorphism
\[
F\Bigl[t,\frac1t,\frac1{1-t}\Bigr]\to F\Bigl[t,\frac1t,\frac1{1-t}\Bigr],  \quad t\mapsto \frac1t,
\]
gives an automorphism $\iota$ of $\dbP^1_F\setminus\{0,1,\infty\}$.
Let $0,1\neq z\in F$, and consider the tangential base point $\overline{0z}=0+z\eps$.
As shown in \cite{Wickelgren12b}*{Lemma 3.16},
\[
(\nu_x)_*(\kappa(\iota(\overline{0z})))=(1/z)_F, \quad (\nu_y)_*(\kappa(\iota(\overline{0z})))=(-1/z)_F.
\]
Therefore the same argument as in Theorem \ref{the main theorem}, combined with Remark \ref{automorphism of U T}, gives rise to $\rho_w^z\in H^1(G,\dbU_{\widetilde T_{w,x}})$ such that in $H^1(G,\calA^+)$,
\[
(\pr_{i,i+1})_*(\rho_w^z)
=\begin{cases}(z)_F,&\hbox{if }a_i=x,\\
(-z)_F,&\hbox{if }a_i=y,
\end{cases}
\qquad i=1,2\nek n.
\]
\end{rem}

\begin{rem}
\rm
The assumption that $n!$ is invertible in $\calA$ can be slightly relaxed, as it is enough to assume that $(j-i)!$ is invertible in $\calA$ for every $(i,j)\in\widetilde T_{w,x}$.
Indeed, these are the denominators that actually appear in the definition of $\bar\varphi_w$.
\end{rem}

\section{Massey Products}
\label{section on Massey Products}
For completeness, we extend the discussion in \cite{Wickelgren12b}, which explains the connection between $1$-cocycles $G\to\dbU$ and Massey products, to the more general setting of the graded Lie algebra $\gr_\bullet(Z,n)$.
First we note the following lemma, whose proof is straightforward:

\begin{lem}
\label{equivalent conditions for 1-cocycle}
Suppose that the profinite group $G$ acts continuously on the profinite group $U$.
A continuous map $\rho\colon G\to U$ is a $1$-cocycle if and only if the map $\rho\rtimes \id\colon G\to U\rtimes G$, $g\mapsto(\rho(g),g)$ is a group homomorphism.
\end{lem}

We refer to \cite{NeukirchSchmidtWingberg}*{Th.\ 1.2.4} for  a detailed description of the Schreier correspondence between equivalence classes of extensions of profinite groups and the second (abelian) cohomology group.
We will need the following basic fact:

\begin{lem}
\label{lemma on Schreier correspondence}
Consider the commutative diagram of profinite groups and continuous maps
\[
\xymatrix{
&&G\ar[d]^{\bar\rho}\ar[ld]_{\rho}&\\
1\to A\ar[r]&B\ar[r]&C\ar[r]&1,
}
\]
where all maps except possibly $\rho$ are homomorphisms, $A$ is abelian, and the row is exact.
Then the pullback to $H^2(G,A)$ of the classifying element in $H^2(C,A)$ of the lower extension is represented by the continuous $2$-cocycle
\[
G\times G\to A, \qquad (g,h)\mapsto \rho(g)\cdot\rho(h)\cdot \rho(gh)\inv.
\]
\end{lem}
\begin{proof}
The pullback of the extension to $G$ along $\bar\rho$ is
\[
1\to A\to B\times_CG\to G\to 1,
\]
where $B\times_CG$ denotes the fiber product, and we view $A$ as a subgroup of $B\times\{1\}$.
The map $g\mapsto (\rho(g),g)$ is a continuous section of the right epimorphism.
Hence this pullback is represented by the continuous $2$-cocycle
\[
(g,h)\mapsto  (\rho(g),g)(\rho(h),h) (\rho(gh),gh)\inv=(\rho(g)\rho(h)\rho(gh)\inv,1).
\qedhere
\]
\end{proof}

Next we consider a subset $Z$ of $\{1,2\nek n\}$ and positive integers $s\leq t$.
By (\ref{decreasing sequence of the V subgroups}) and Proposition \ref{commutativity of product}, the group $V_{T(Z,t)}/V_{T(Z,t+1)}$ is abelian, and there is a central extension of profinite groups
\begin{equation}
\label{central extension}
1\to V_{T(Z,t)}/V_{T(Z,t+1)}\to V_{T(Z,s)}/V_{T(Z,t+1)}\to V_{T(Z,s)}/V_{T(Z,t)}\to 1.
\end{equation}
Note that $V_{T(Z,t)}/V_{T(Z,t+1)}\isom \Ker(\dbU_{T(Z,t+1)}\to\dbU_{T(Z,t)})$, etc.

Let $G$ be as before, a profinite group which acts on the pro-$\calC$ ring $\calA$ via a continuous homomorphism $\chi\colon G\to\hat\dbZ(\calC)^\times$, and let $G$ act on $\dbU=\dbU_n(\calA)$, and therefore on its quotients $\dbU_T$, as in \S\ref{section on equivariant maps}.

\begin{prop}
\label{lifting of 1-cocycle}
Let $\rho\colon G\to  V_{T(Z,s)}/V_{T(Z,t+1)}$ be a continuous map such that the induced map $\bar\rho\colon G\to  V_{T(Z,s)}/V_{T(Z,t)}$ is a $1$-cocycle.
Then
\begin{enumerate}
\item[(a)]
The map
\[
c\colon G\times G\to V_{T(Z,t)}/V_{T(Z,t+1)}, \quad (g,h)\mapsto \rho(g)\cdot{}^g\rho(h)\cdot\rho(gh)\inv
\]
is a well-defined continuous $2$-cocycle.
\item[(b)]
There exists a continuous $1$-cocycle $G\to V_{T(Z,s)}/V_{T(Z,t+1)}$ lifting $\bar\rho$ if and only if $c$ is cohomologous to $0$.
\end{enumerate}
\end{prop}
\begin{proof}
From (\ref{central extension}) and Lemma \ref{equivalent conditions for 1-cocycle} we obtain a central extension of profinite groups as in the following commutative diagram:
\begin{equation}
\label{embedding problem}
\xymatrix{&&&G\ar[d]^{\bar\rho\rtimes\id}\ar[dl]_{\rho\rtimes\id}&\\
1\ar[r]& \dfrac{V_{T(Z,t)}}{V_{T(Z,t+1)}}\ar[r]& \dfrac{V_{T(Z,s)}}{V_{T(Z,t+1)}}\rtimes G\ar[r]&  \dfrac{V_{T(Z,s)}}{V_{T(Z,t)}}\rtimes G\ar[r]& 1,
}
\end{equation}
where all maps except possibly $\rho\rtimes\id$ are homomorphisms.

Since $\bar\rho$ is a $1$-cocycle, $\rho(g)\cdot{}^g\rho(h)\cdot\rho(gh)\inv$ is in the kernel $V_{T(Z,t)}/V_{T(Z,t+1)}$ of  the projection $ V_{T(Z,s)}/V_{T(Z,t+1)}\to V_{T(Z,s)}/V_{T(Z,t)}$ for every $g,h\in G$, showing that $c$ is well defined.

By Lemma \ref{lemma on Schreier correspondence}, the pullback to $G$ of the central extension of (\ref{embedding problem}) along the homomorphism $\bar\rho\rtimes\id$ is represented by the continuous $2$-cocycle
\[
(g,h)\mapsto (\rho\rtimes\id)(g)\cdot(\rho\rtimes\id)(h)\cdot(\rho\rtimes\id)(gh)\inv
=(\rho(g)\cdot{}^g\rho(h)\cdot\rho(gh)\inv,1),
\]
which we may identify with $c$.
In particular, this gives (a).

By Lemma \ref{equivalent conditions for 1-cocycle}, there is a $1$-cocycle as in (b), if and only if there is $\rho$ as above, such that $\rho\rtimes\id$ is a homomorphism, i.e., the embedding problem (\ref{embedding problem}) is solvable.
By Hoechsmann's lemma  \cite{NeukirchSchmidtWingberg}*{Prop.\ 3.5.9}, this means that  the pullback to $G$ of the central extension is a trivial extension.
By what we have seen, this means that $[c]=0$.
\end{proof}

We call the cohomology class
\[
[c]\in H^2(G, V_{T(Z,t)}/V_{T(Z,t+1)})
\]
the \textsl{generalized Massey product of $Z$ at levels $s\leq t$ and corresponding to the $1$-cocycle $\bar \rho$}.
This terminology is motivated by the following example:

\begin{exam}
\label{Example Z 12...n}
\rm
Take $Z=\{1,2\nek n\}$, $s=1$ and $t=n$.
We have
\[
T(Z,1)=\Del,  \quad  T(Z,n)=\dbI\setminus\{(1,n+1)\}, \quad T(Z,n+1)=\dbI.
\]
Hence (\ref{central extension}) becomes
\[
0\to \calA\to\dbU\to\overline{\dbU}\to 1,
\]
where $\overline{\dbU}=\dbU_{\dbI\setminus\{(1,n+1)\}}$, i.e., it is $\dbU$ with the upper-right entry omitted.
Then $\calA$ embeds in $\dbU$ via $\alp\mapsto I+\alp E_{1,n+1}$, where $E_{1,n+1}$ is the matrix which is $1$ at entry $(1,n+1)$ and is $0$ elsewhere.
Here the $G$-action on $\calA$ is $n$-twisted, i.e., is given by ${}^g\alp=\chi(g)^n\cdot\alp$.

For $\rho$, $\bar\rho$, and $c$ as in Proposition \ref{lifting of 1-cocycle}, and for $g,h\in G$ we get
\[
(I+c(g,h)E_{1,n+1})\rho(gh)=\rho(g)\cdot{}^g\rho(h).
\]
At entry $(1,n+1)$ this gives
\[
\rho(gh)_{1,n+1}+c(g,h)=\sum_{k=1}^{n+1}\rho(g)_{1,k}\cdot\chi(g)^{n+1-k}\rho(h)_{k,n+1}.
\]
Writing $\partial\colon C^1(G,\calA)\to C^2(G,\calA)$ for the coboundary map on continuous cochains, the latter equality becomes
\[
c(g,h)=\sum_{k=2}^n(\rho_{1k}\cup\rho_{k,n+1})(g,h)+(\partial \rho_{1,n+1})(g,h).
\]
It follows that $\sum_{k=2}^n\rho_{1k}\cup\rho_{k,n+1}$ is also a continuous $2$-cocycle, which is cohomologous to $c$.

We conclude from Proposition \ref{lifting of 1-cocycle} that the continuous $1$-cocycle $\bar\rho\colon G\to\overline\dbU_n$ lifts to a continuous $1$-cocycle $G\to \dbU_n$ if and only if
\[
\Bigl[\sum_{k=2}^n\rho_{1k}\cup\rho_{k,n+1}\Bigr]=0\in H^2(G,\calA).
\]

The left-hand side here is the element of the \textsl{$n$-fold Massey product} corresponding to $\bar\rho$, considered as a \textsl{defining system}, as studied in \cite{Wickelgren12b}, \cite{Wickelgren 17}.
When the $G$-action on $\calA$ is trivial, this analysis goes back to Dwyer \cite{Dwyer75} (in the discrete group case; see e.g., \cite{Efrat14} for the profinite case).
\end{exam}

\begin{exam}
\label{example n 2}
\rm
Take in the previous example $n=2$.
Then a continuous $1$-cocycle $\bar\rho\colon G\to\overline\dbU_2$ consists of the $1$-cocycles $\bar\rho_{12},\bar\rho_{23}$.
We obtain that it lifts to a continuous $1$-cocycle $G\to\dbU_2$ if and only if $[\bar\rho_{12}]\cup[\bar\rho_{23}]=0$ in $H^2(G,\calA)$, where $G$ acts on $\calA$ by ${}^g\alp=\chi(g)^2\alp$.
\end{exam}

\end{document}